\documentclass[reqno,11pt]{article}
\usepackage{amsmath, latexsym, amsfonts, amssymb, amsthm, amscd}
\usepackage{graphics,epsf,psfrag,enumitem}
\usepackage{graphics,epsf,psfrag}
\usepackage[utf8]{inputenc}
\usepackage{stmaryrd}
\usepackage{hyperref}
\usepackage{xcolor}

\usepackage{authblk}

\usepackage{enumitem}
\setlist{itemsep=1pt,topsep=2pt,parsep=1pt}

\numberwithin{equation}{section}

\let\OLDthebibliography\thebibliography
\renewcommand\thebibliography[1]{
  \small
  \OLDthebibliography{#1}
  \setlength{\parskip}{0pt}
  \setlength{\itemsep}{3pt plus 0.3ex}
}

\newtheorem{theorema}{Theorem}

\newtheorem{theorem}{Theorem}[section]

\newtheorem{proposition}[theorem]{Proposition}

\newtheorem{remark}[theorem]{Remark}
\newtheorem{example}[theorem]{Example}

\newcommand{\ind}{\mathbf{1}}

\renewcommand{\hat}{\widehat}
\newcommand{\dd}{\mathrm{d}}

\newcommand{\bP}{{\ensuremath{\mathbf P}} }

\newcommand{\bE}{{\ensuremath{\mathbf E}} }

\newcommand{\bbN}{{\ensuremath{\mathbb N}} }

\newcommand{\bbR}{{\ensuremath{\mathbb R}} }

\newcommand{\gep}{\varepsilon}       

\newcommand{\go}{\omega}

\newcommand{\hookdoubleheadrightarrow}{%
  \hookrightarrow\mathrel{\mspace{-15mu}}\rightarrow
}

\makeatletter
\def\namedlabel#1#2{\begingroup
   \def\@currentlabel{#2}%
   \label{#1}\endgroup
}
\makeatother

\begin{document}

\title{An application of Sparre Andersen's fluctuation theorem\\
 for exchangeable and sign-invariant random variables}


\author[*,$\dagger$]{Quentin Berger}
\author[*]{Lo\"ic B\'ethencourt}
\affil[*]{\footnotesize Laboratoire de Probabilit\'e Statistique et Mod\'elisation, Sorbonne Universit\'e.}
\affil[$\dagger$]{\footnotesize D\'epartement de Math\'ematiques et Applications, \'Ecole Normale Sup\'erieure, PSL.}

\date{}

\maketitle

\begin{abstract}
\noindent
We revisit the celebrated Sparre Andersen's fluctuation result on persistence (or survival) probabilities $\bP(S_k> 0 \;\forall\, 0\leq k\leq n)$ for symmetric random walks $(S_n)_{n\geq 0}$. We give a short proof of this result when  considering sums of random variables that are only assumed \emph{exchangeable} and \emph{sign-invariant}. We then apply this result to the study of persistence probabilities of (symmetric) additive functionals of Markov chains, which can be seen as a natural generalization of integrated random walks.\\[3pt]
\textit{MSC2020 AMS classification:} 60G50; 60J55; 60G09; 60J10.
\end{abstract}

\section{Introduction}

We discuss in this paper several results on persistence problems for symmetric random walks, in particular the well known (but nonetheless striking) fluctuation theorem due to Sparre Andersen~\cite{Sparre54}.
To introduce the result, let $\xi_1, \ldots, \xi_n$ be real random variables and let us denote $S_0=0$, and 
$S_k=\sum_{i=1}^k \xi_i$ for $1\leq k \leq n$.
Then, one is interested in estimating the persistence probabilities
\[
\bP\big(S_k > 0 \, \text{ for all } 1\leq k \leq n \big) \quad\text{ or } \quad \bP\big(S_k \geq 0 \, \text{ for all } 1\leq k \leq n \big)\,.
\]
Sparre Andersen's result~\cite{Sparre54} states that if the $(\xi_i)_{1\leq i\leq n}$ are i.i.d., symmetric and have no atom, then these probabilities do not depend on the law of $\xi$.

\begin{theorema}
\label{thmA}
If the $(\xi_i)_{1\leq i \leq n}$ are i.i.d.\ and if the distribution of $\xi_i$ is symmetric and has no atom, then we have 
\[
\bP\big(S_k > 0 \, \text{ for all } 1\leq k \leq n \big) = g(n) =\bP\big(S_k \geq 0 \, \text{ for all } 1\leq k \leq n \big) \,,
\]
where
\[
g(n) := \frac{(2n-1)!!}{(2n)!!}  = \frac{1}{4^n} \binom{2n}{n} =  \prod_{k=1}^n \Big( 1- \frac{1}{2k} \Big) \,.
\]
\end{theorema}

In fact, \cite{Sparre54} proves a number of results for exchangeable random variables, but the independence plays an important role in obtaining fluctuation's results, such as arcsine laws or persistence probabilities.
Chapter XII in Feller's book~\cite{FellerII} provides a streamlined approach to Sparre Andersen's result, using a duality argument together with a (purely combinatorial) cyclic lemma, see~\cite[XII.6]{FellerII}: these two ingredients only require  the law of $(X_1,\ldots,X_n)$ to be exchangeable.
Independence is required in the last step of the proof, in order to obtain that the (weak) ladder epochs defined iteratively by $T_0=0$ and $T_{k} = \min\{n > T_{k-1}, S_n \geq  \max_{0\leq j<n}\{S_j\} \}$ 
form a renewal sequence.
The following is then deduced, see~\cite[XII.7, Thm.~1]{FellerII}:
\begin{theorema}
\label{thm:B}
If the $(\xi_i)_{i\geq 1}$ are i.i.d., then for any $s\in [0,1)$ we have
\[
1- \bE\big[s^{T_1}\big] = \exp\big( - \sum_{n=1}^{\infty} \frac{s^n}{n} \bP(S_n \geq 0) \big) \,.
\]
\end{theorema}

\noindent
Theorem~\ref{thm:B} can be seen as a particular case of a Wiener--Hopf factorization, also known as Spitzer--Baxter formula, which gives the joint Laplace transform/characteristic function of the first ladder epoch and ladder height, see~\cite[XII.9]{FellerII}.
Note that if the law of $X_i$ is symmetric and has no atom, one has $\bP(S_n \geq 0)=\frac12$ for all $n$, so the generating function of $T_1$ is given by 
$\bE[s^{T_1}] =  1- \sqrt{1-s}$.

\smallskip
Let us mention that another proof of Theorem~\ref{thmA} is presented in~\cite[Prop.~1.3]{DDG13}: the proof is remarkably simple and elegant and relies on a decomposition with respect to the first time that $(S_k)_{k \geq 0}$ hits its minimum $\min_{1\leq j\leq n}\{S_j\}$, using the independence of the $(\xi_i)_{1\leq i \leq n}$.
The goal of our paper is to present a version of Theorem~\ref{thmA} valid for exchangeable and sign-invariant random variables $(\xi_i)_{1\leq i \leq n}$, going beyond the independent and symmetric setting. Sparre Andersen was already aware of this result, see \cite[Thm.~4]{Sparre54} (its proof is however a bit laconic), but we give here a short and self-contained  proof (with no combinatorial lemma), taking inspiration from~\cite{DDG13}.
We then present some application to persistence probabilities for symmetric additive functionals of birth-death chains, giving in particular a simple proof of Sinai's result~\cite{Sinai92} on the integrated simple random walk.

\section{The case of exchangeable and sign-invariant random variables}

Below, we consider examples in which there is no independence: one is however still able to obtain the same conclusion as in Theorem~\ref{thmA}, using some weaker but natural assumption. Let us introduce some terminology: the law $\bP$ of $(\xi_1,\ldots, \xi_n)$ is said to be:
\begin{description}
\item[(E)] \namedlabel{hyp:E}{(E)}
exchangeable if for any permutation $\sigma \in \mathfrak{S}_n$,  $(\xi_{\sigma(1)},\ldots, \xi_{\sigma(n)})$ has the same distribution as $(\xi_1,\ldots, \xi_n)$;
\item[(S)] \namedlabel{hyp:S}{(S)}
sign-invariant if for any $\gep = (\gep_1,\ldots, \gep_n)\in \{-1,1\}^n$,
$(\gep_1 \xi_{1},\ldots, \gep_n \xi_{n})$ has the same distribution as $(\xi_1,\ldots, \xi_n)$.
\end{description}
One can think about these two assumptions as an invariance of $\bP$ under the action of the groups: of permutations $\mathfrak{S}_n$ for condition~\ref{hyp:E}; of sign changes $\{-1,1\}^{n}$ for condition~\ref{hyp:S}.

The following result is actually already mentioned in \cite[Thm.~4]{Sparre54}: the present paper is meant as a way to put more emphasis on this result;
we then provide some applications.
\begin{theorem}
\label{thm:exchang}
Assume that $(\xi_1, \ldots, \xi_n)$ satisfies conditions~\ref{hyp:E}-\ref{hyp:S}.
Then, recalling that we have defined $g(n)= \frac{1}{4^n} \binom{2n}{n}$, we have
\[
\bP\big(S_k > 0 \, \text{ for all } 1\leq k \leq n \big) \leq g(n) \leq \bP\big(S_k \geq  0 \, \text{ for all } 1\leq k \leq n \big) \,.
\]
Notably, if $\bP$ has no atom,
$\bP\big(S_k > 0  \text{ for all } 1\leq k \leq n \big)  =  g(n) = \bP\big(S_k \geq 0  \text{ for all } 1\leq k \leq n \big)$.
\end{theorem}

\begin{remark}
\label{rem:bounds1}
Let us note that one may also obtain a lower bound  on $\bP(S_k > 0 \, \text{ for all } 1\leq k \leq n )$ by imposing that $S_1>0, S_2-S_1\geq 0, \ldots, S_{n}-S_1\geq 0$: since conditionally on $X_1$ the random vector $(X_2,\ldots, X_n)$ is still exchangeable and sign-invariant, we get from Theorem~\ref{thm:exchang}
\begin{equation}
\label{eq:bound}
\bP(S_k > 0 \, \text{ for all } 1\leq k \leq n ) \geq \bP(X_1 >0) g(n-1) \,,
\end{equation}
so we have both an upper and a lower bound on $\bP(S_k > 0 \, \text{ for all } 1\leq k \leq n )$ in terms of~$g(\cdot)$.
Let us also stress that we have $ \frac{1}{\sqrt{\pi (n+1/2)}}\leq g(n) \leq \frac{1}{\sqrt{\pi n}}$ for all $n\geq 1$, see e.g.~\cite[Eq.~(3.1)]{BDGJS23}.

\noindent
As far as $\bP(S_k \geq 0 \text{ for all } 1\leq k \leq n )$ is concerned, Theorem~\ref{thm:exchang} provides a lower bound, but there is no general upper bound: \cite[Conj.~6]{BDGJS23} conjectures that this probability is maximal when $(S_k)_{k\geq 0}$ is the simple symmetric random walk, \textit{i.e.}\ that $\bP(S_k \geq 0 \text{ for all } 1\leq k \leq n ) \leq g( \lfloor n/2 \rfloor)$.
\end{remark}

\noindent
One can easily see that both conditions~\ref{hyp:E}-\ref{hyp:S} are necessary to obtain the statement of Theorem~\ref{thm:exchang}, see for instance Section~\ref{sec:counterexample} below.

\begin{remark}[About sign-invariance]
The sign-invariance of $(\xi_1, \ldots, \xi_n)$ can be seen as a strong form of symmetry for the law of $(\xi_1, \ldots, \xi_n)$.
Let us mention the article \cite{Ber65} which derives some properties of sign-invariant sequences: in particular, \cite[Lem.~1.2]{Ber65} tells that sign-invariant $(\xi_1, \ldots, \xi_n)$ are independent conditionally on $(|\xi_1|, \ldots, |\xi_n|)$.
In other words, any exchangeable and sign-invariant law can be obtained by taking a random vector $(Z_1,\ldots, Z_n)$ (that can be constrained to have $0\leq Z_1\leq \cdots \leq Z_n$) and by shuffling and changing the signs of its coordinates randomly, that is setting $\xi_i = \gep_i Z_{\sigma(i)}$ for $1\leq i \leq n$, where the $(\gep_i)_{1\leq i \leq n}$ are i.i.d.\ signs (\textit{i.e.}\ such that $\bP(\gep_i=\pm1) =\frac12$) and $\sigma$ is a random uniform permutation of $\{1,\ldots, n\}$ (\textit{i.e.}\ such that $\bP(\sigma =\nu) = \frac{1}{n!}$ for all $\nu\in \mathfrak{S}_n$). 
\end{remark}

\begin{remark}[Finite vs.\ infinite sequences]
The law of an infinite sequence  $(\xi_i)_{i\geq 1}$ of random variables is called exchangeable, resp.\ sign-invariant, if for any $n\geq 1$ the law of $(\xi_1,\ldots, \xi_n)$ is exchangeable, resp.\ sign-invariant.
By de Finetti's theorem, one knows that the law of $\xi=(\xi_i)_{i\geq 1}$ is exchangeable if and only if it is conditionally i.i.d.: in other words,
there exists a \textit{random} probability distribution $\mu$ such that $\bP(\xi \in \cdot \mid \mu) = \mu^{\otimes\infty}$, where $\mu^{\otimes \infty}$ denotes the law of an infinite sequence of i.i.d.\ random variables with law $\mu$. 
In~\cite[Thm.~1]{Ber62}, it is shown that an exchangeable sequence $(\xi_i)_{i\geq 1}$ is sign-invariant if and only if $\mu$ is almost surely symmetric.
From this, it is easy to derive Theorem~\ref{thm:exchang} for an exchangeable and sign-invariant sequence $(\xi_i)_{i\geq 1}$, simply by working conditionally on the realization of $\mu$. Hence, the truly remarkable fact about Theorem~\ref{thm:exchang} is that it holds for a \textit{finite} exchangeable and sign-invariant $(\xi_1,\ldots, \xi_n)$.
\end{remark}

\section{Persistence probabilities of $f$-integrated birth-death chains}
\label{sec:integrated}

Let us consider $(X_n)_{n\geq 0}$ a birth and death Markov chain, starting from $X_0=0$: it is a Markov chain on $\mathbb Z$ such that $|X_n-X_{n-1}|\leq 1$, with transition probabilities
\[
p_x = p(x,x+1),\ \  q_x= p(x,x-1),\ \  r_x = p(x,x) = 1-p_x-q_x \quad \text{ for } x\geq 1.
\]
We assume that $(X_n)_{n\geq 0}$ is symmetric, in the sense that $(-X_n)_{n\geq 0}$ has the same distribution as $(X_n)_{n\geq 0}$; in other words, the transition probabilities verify $p(x,y) = p(-x,-y)$ for any $x,y\in \mathbb Z$.
For $x=0$, we also set $p_0 = p(0,1) \in (0,\frac12]$, $q_0=p(0,-1) =p_0$ and $r_0=p(0,0)=1-2p_0$.

Let $f:\mathbb Z \to \mathbb R$ be any anti-symmetric function which preserves the sign of $x$, \textit{i.e.}\ such that $f(x)> 0$ for $x> 0$ and $f(x)< 0$ for $x< 0$ (and naturally $f(0)=0$).
Then we define the \textit{$f$-integrated} Markov chain, or \emph{additive functional}, as follows: $\zeta_0=0$ and, for $n\geq 1$,
\begin{equation}
\label{def:zeta}
\zeta_n := \sum_{i=1}^n f(X_i) \,.
\end{equation}
We are now interested in the persistence (or survival) probabilities 
\begin{equation*}
\bP(\zeta_k \geq 0 \text{ for all } 1\leq k \leq n)
\quad \text{ or } \quad
\bP(\zeta_k > 0\text{ for all } 1\leq k \leq n)\,.
\end{equation*}
A classical, well-studied example is when $(X_n)_{n\geq 0}$ is the simple symmetric random walk and $f$ is the identity: then $(\zeta_n)_{n\geq 0}$ is the integrated random walk and the persistence probabilities are known to be of order $n^{-1/4}$, see~\cite{Sinai92}.

As another motivation, let us point out that the persistence problem for integrated random walk bridges appeared in the context of a polymer pinning model \cite{CD09} (see also~\cite{ADL14}) and more recently in the study of graphic sequences in~\cite{BDGJS23}, \textit{i.e.}\ on the number $G(n)$ of integer sequences $n-1\geq d_1\geq \cdots \geq d_n \geq 0$ that are the degree sequence of a graph.
We are therefore also interested in the behavior of
\begin{equation*}
\bP(\zeta_k > 0 \text{ for all } 1\leq k \leq n, X_{n}=0)
\quad \text{ or } \quad
\bP(\zeta_k \geq 0\text{ for all } 1\leq k \leq n, X_{n}=0)\,.
\end{equation*}
(Let us assume that $(X_n)_{n\geq 0}$ is aperiodic for the simplicity of exposition.)

\begin{theorem}
\label{thm:integrated}
Let $(X_n)_{n\geq 0}$ be a symmetric recurrent birth-death chain and let $f$ be an odd function such that $xf(x)\geq 0$; recall that $p_0:=\bP(X_1>0)$.
Then, recalling that $g(n):=\frac{1}{4^n} \binom{2n}{n}$ and setting by convention $g(n)=0$ for $n\leq 0$, we have for all $n\geq 1$
\[
 p_0(1-p_0) \bE\big[  g(L_n-1) \big] \leq \bP\big(\zeta_k > 0 \text{ for all } 1\leq k \leq n \big)\leq \bE\big[  g(L_n) \big]
\]
where $L_n := \sum_{k=1}^n \ind_{\{X_k=0\}}$ is the local time of the chain at $0$ up to time $n$. Regarding the bridge, we have for every $n\geq1$
\[
 p_0 \bE\big[  g(L_n-1) \ind_{\{X_n=0\}} \big]  \leq \bP\big(\zeta_k > 0 \text{ for all } 1\leq k \leq n , X_n =0\big) \leq \bE\big[  g(L_n) \ind_{\{X_n=0\}} \big].
\]
\end{theorem}

\noindent
Let us mention that our approach is based on an excursion decomposition of the process $(X_n)_{n\geq 0}$.
Our results would still be valid for Markov chains that ``cannot jump above $0$'', \textit{i.e.}\ such that $\bP(X_n\leq 0 \mid X_{n-1}=x) = \bP(X_n=0\mid X_{n-1}=x)$ for any $x\in \mathbb N$.
We have chosen to stick with the birth-death setting since it is the most natural example of such chains (and already contains a wide class of behaviors).

Since we have the asymptotic behavior $g(m)\sim (\pi m)^{-1/2}$ as $m\to\infty$ (we actually have explicit bounds $ (\pi (m+1/2))^{-1/2} \leq g(m) \leq (\pi m)^{-1/2}$ for $m\geq 1$, see Remark~\ref{rem:bounds1}), we can give sufficient conditions, often verified in practice, to obtain the asymptotic behavior of $\bE[g(L_n)]$, $\bE[g(L_n)\ind_{\{X_n=0\}}]$ as $n\to\infty$.
(The same asymptotics hold with $g(L_n-1)$ in place of $g(L_n)$.)


\begin{proposition}
\label{lem:Ln}
Assume that $(X_n)_{n\geq 0}$ is aperiodic and let  $\tau_1 := \inf\{n\geq 1, X_n=0\}$.

(i) If $(X_n)_{n\geq 0}$ is positive recurrent, then we have, as $n\to\infty$
\begin{equation}
\label{eq:persist1}
\bE[g(L_n)] \sim \frac{\sqrt{\bE[\tau_1]}}{\sqrt{\pi n}} \qquad \text{ and } \qquad \bE[g(L_n)\ind_{\{X_{n}=0\}}] \sim \frac{1}{\sqrt{\pi n \bE[\tau_1]}} \,.
\end{equation}

(ii) If $\bP(\tau_1>n) \sim \ell(n) n^{-\alpha}$ for some $\alpha \in [0,1]$ and some slowly varying function $\ell(\cdot)$, then we have, as $n\to\infty$
\begin{equation}
\label{eq:persist2}
\bE[g(L_n)] \sim  \frac{c_{\alpha}}{\sqrt{\pi b_n}} \quad 
\text{ with } \  b_n := 
\begin{cases}
n/\mu(n) \\
n^{\alpha}/\ell(n) \\
1/\ell(n)
\end{cases}
 c_{\alpha} := 
\begin{cases}
1 &  \text{ if } \alpha=1 \,,\\
\Gamma(1-\alpha)^{1/2}\bE\big[(\mathcal Z_{\alpha})^{\alpha/2} \big] & \text{ if } \alpha\in (0,1) \,,\\
\sqrt{\pi} & \text{ if } \alpha=0\,,
\end{cases}
\end{equation}
where $\mu(n):=\bE[\tau_1 \ind_{\{\tau_1\leq n\}}]$ and with $\mathcal Z_{\alpha}$ a one-sided $\alpha$-stable random variable of Laplace transform~$e^{- t^\alpha}$.
Moreover, if $\alpha \in (2/3,1)$, we have the asymptotic
\begin{equation}
\label{eq:persistlocal}
\bE[g(L_n)\ind_{\{X_{n}=0\}}] \sim c'_{\alpha} \, n^{\frac{\alpha}{2} -1} \ell(n)^{-1/2}\,, \qquad \text{ with } \ c'_{\alpha} = \alpha\Gamma(1-\alpha)^{-1/2}\bE\big[(\mathcal Z_{\alpha})^{-\alpha/2} \big]\,.
\end{equation}
In order to have the same asymptotic~\eqref{eq:persistlocal} for $\alpha\in(0, 2/3]$, a sufficient condition is that there exists a constant $C>0$ such that $\bP(\tau_1 = n) \leq C \ell(n) n^{-(1+\alpha)}$ for all $n\geq 1$.
\end{proposition}

\begin{example}
We give below a class of example where Proposition~\ref{lem:Ln} can be applied, but let us comment on the case of the simple symmetric random walk $(X_n)_{n\geq 0}$.
In that case, Proposition~\ref{lem:Ln} is verified (up to periodicity issues).
We then get that the persistence probabilities for the integrated random walk verify
\begin{equation}
\label{eq:boundRW}
 (1+o(1)) \tfrac{1}{4} c_0 n^{-1/4}\leq \bP\big(\zeta_k >0 \text{ for all } 1\leq k \leq 2n \big) \leq  (1+o(1)) c_0 n^{-1/4} \,,
\end{equation}
with $c_0$ an explicit constant. This recovers a result by Sinai~\cite{Sinai92} (in the case $f(x)=x$).
For the integrated simple random walk bridge, we also have that
\begin{equation}
\label{eq:boundRW}
 (1+o(1)) \tfrac{1}{4} c_1 n^{-3/4}\leq \bP\big(\zeta_k >0 \text{ for all } 1\leq k \leq 2n , X_{2n}=0\big) \leq  (1+o(1)) c_1 n^{-3/4} \,,
\end{equation}
recovering a result by Vysotsky~\cite{Vysotsky14}.
\end{example}

\begin{remark}
In the case $\alpha \in (0,1)$ the condition $\bP(\tau_1>n)\sim \ell(n) n^{-\alpha}$ is exactly equivalent to $\tau_1$ being in the domain of attraction of an $\alpha$-stable law.
Our sufficient condition $\alpha>2/3$ to get~\eqref{eq:persistlocal} may seem mysterious at first, but it is reminiscent of the Garsia--Lamperti~\cite{GL62} condition $\alpha >1/2$ to obtain a Strong Renewal Theorem (SRT), \textit{i.e.}\ the sharp asymptotic behavior of $\bP(X_n=0) = \bP(n\in \tau)$.
In the case $\alpha \leq 2/3$, our condition on the local tail $\bP(\tau_1 =n)$ is reminiscent of Doney's condition in~\cite[Eq.~(1.9)]{D97} to obtain a SRT,
Let us also mention~\cite{CD19} where a necessary and sufficient condition for a SRT is found: with some effort, it should translate in a necessary and sufficient condition for~\eqref{eq:persistlocal} to hold.
Note that for technical simplicity we only deal with the case $\alpha\in (0,1)$ in~\eqref{eq:persistlocal}, but an analogous result should hold in the cases $\alpha=1$ and $\alpha=0$.
\end{remark}

\begin{remark}
We stress here that the asymptotic bounds that we obtain combining Theorem~\ref{thm:integrated} and Proposition~\ref{lem:Ln} do not depend on the (anti-symmetric) function $f$ in the definition~\eqref{def:zeta} of $\zeta_n$. 
In particular, the bounds~\eqref{eq:boundRW} obtained in the case of the simple random walk (and for Bessel-like random walks, see~\eqref{eq:persistBRW1}-\eqref{eq:persistBRW2} below) are valid for \textit{any} odd function function $f$. 
\end{remark}

\subsection*{A class of examples: Bessel-like random walks}

Let us give more precise result in the case of symmetric Bessel-like random walk, see~\cite{Ale11,HR53,Lamp62}.
This is a class of birth-death chains that includes the simple symmetric random walk. They have the following transition probabilities: for $x\geq 1$
\begin{equation}
\label{def:bessel}
p_x:=p(x,x+1) = \frac12 \Big(1-\frac{\delta+\gep_x}{2x}\Big) \,, \qquad p(x,x-1) = 1-p_x \,,
\end{equation}
where $\delta\in \mathbb R$ and $\gep_x$ is such that $\lim_{x\to\infty} \gep_x =0$;
we take $p(-x,-x-1) = p(x,x+1)$ and $p(0,1)=p(0,-1)=\frac12$, for symmetry reasons.
We also assume uniform ellipticity, \textit{i.e.}\ there is some $\eta>0$ such that $p_x\in [\eta,1-\eta]$ for all $x\in \mathbb Z$.
The parameter~$\delta$ is called the \textit{drift parameter} and we have the following behavior: the walk $(X_n)_{n\geq 0}$ is transient if $\delta<-1$, recurrent if $\delta>-1$, positive recurrent if $\delta>1$; in the cases $\delta=-1$ and $\delta =1$, the behavior depends on the function~$\gep_x$.

More precisely, letting $\lambda_x = \prod_{k=1}^x \frac{p_x}{1-p_x}$, the random walk $(X_n)_{n\geq 1}$ is recurrent if and only if $\sum_{x=1}^{\infty} \lambda_x=+\infty$.
Moreover, with the notation~\eqref{def:bessel}, there is a constant $K_0>0$ such that
\begin{equation*}
\label{eq:lambda}
\lambda_x = \prod_{k=1}^x \frac{1-p_k}{p_k} \sim K_0 \,x^{\delta}L(x)^{-1} \quad \text{ as } x\to\infty,\quad  \text{ where } L(x) := \exp\Big(\sum_{k=1}^x \frac{\gep_k}{k} \Big)\,.
\end{equation*}
Note that $L$ is a slowly varying function.
Then, in~\cite[Thm.~2.1]{Ale11}, the sharp tail of the distribution of $\tau_1$ is derived: 
Let $\delta \geq -1$ and set $\alpha:=\frac{1+\delta}{2}$, then as $n\to\infty$
\begin{equation*}
\label{eq:alexander}
\begin{split}
\text{ if } \delta>-1\ (\alpha>0), &\qquad  \bP(\tau_1 > n) \sim \tfrac{2^{1-\alpha}}{K_0\Gamma(\alpha)}\; n^{-\alpha} L(\sqrt{n})  \,,\\
\text{ if } \delta=-1\ (\alpha=0), &\qquad  \bP(\tau_1 > n) \sim \tfrac{1}{K_0}\; \nu(n)^{-1} \quad \text{where } \nu(n):=\sum\limits_{x\leq n,\, x\text{ even}} \frac{1}{x L(\sqrt{x})} \,.
\end{split}
\end{equation*}
Also, if $\delta=-1$, $(X_n)_{n\geq 0}$ is recurrent if and only if $\lim_{n\to\infty} \nu(n)=\infty$.
 
For $\alpha>0$, \cite[Thm.~2.1]{Ale11} also gives that 
\[
\bP(\tau_1 =n) \sim \tfrac{2^{2-\alpha}}{K_0\Gamma(\alpha)}\; n^{-(1+\alpha)} L(\sqrt{n}) 
\text{ as $n\to\infty$, $n$ even}.
\]

One can then simply apply Proposition~\ref{lem:Ln}  (up to periodicity issues) with $\ell(n) = c L(\sqrt{n})$ or $\ell(n)=c/\nu(n)$ to obtain that, in the null-recurrent case:
\begin{equation}
\label{eq:persistBRW1}
(1+o(1)) \,\tfrac{1}{4}  c_\alpha (b_n)^{-1/2} \leq \bP\big(\zeta_k >0 \text{ for all } 1\leq k \leq 2n \big) \leq  (1+o(1)) \, c_\alpha (b_n)^{-1/2} \,,
\end{equation}
with $b_n$ given in~\eqref{eq:persist2}.
Also, for $\alpha\in (0,1)$,
\begin{equation}
\label{eq:persistBRW2}
  c' n^{\frac{\alpha}{2}-1} L(\sqrt{n})^{-1/2} \leq \bP\big(\zeta_k >0 \text{ for all } 1\leq k \leq 2n , X_{2n}=0\big) \leq  c'' n^{\frac{\alpha}{2}-1} L(\sqrt{n})^{-1/2} \,.
\end{equation}

%
%
%
%

\subsection*{Further comments and comparison with the literature}

Let us stress that there are several directions in which one could extend our results.
First, one could consider more general underlying random walks or Markov chains, for instance random walks with non-necessarily symmetric increments or Markov chains that can ``jump over $0$''.
One could also consider more general functions~$f$, not necessarily symmetric.
Another room for improvement is to obtain sharp asymptotics for the persistence probabilities, \textit{i.e.}\ for instance finding the correct constant $\hat c_{\alpha}$ such that 
$\bP(\zeta_k >0 \text{ for all } 1\leq k \leq n )\sim \hat c_{\alpha} b_{n}^{-1/2}$ in~\eqref{eq:persistBRW1}.

\smallskip
A big part of the literature has considered the case of integrated random walks, that is considering the Markov chain $X_n=\sum_{k=1}^n \xi_k$ with $(\xi_k)_{k\geq 0}$ i.i.d.\ random variables and $\zeta_n=\sum_{i=1}^n X_i$
(\textit{i.e.}\ taking $f(x)=x$), starting with the work of Sinai~\cite{Sinai92}.
Under the condition that the $\xi_k$'s are centered with a finite second moment, the persistence probability $\bP(\zeta_k >0 \text{ for all } 1\leq k \leq n )$ has been proven to be of order $n^{-1/4}$ in~\cite{DDG13}, and the sharp asymptotic  $\sim c_1 n^{-1/4}$ in~\cite{DW15}.
The case where $\xi_k$ does not have a finite second moment remains mostly open, except in some specific (one-sided) cases, see e.g.~\cite{DDG13,Vysotsky14}.
The persistence probabilities of integrated random walk bridges has been studied for instance in~\cite{ADL14,BDGJS23, CD09,Vysotsky14}: the order $n^{-3/4}$ is found in~\cite[Prop.~1]{Vysotsky14} for centered random walks with finite variance and the asymptotic $\sim c n^{-3/4}$ is given in~\cite[Prop.~1.2]{BDGJS23} in the case of the simple random walk.
Our result can be seen as an extension to persistence problems for additive functionals of  Markov chains $(X_n)_{n\geq 0}$; with the major restriction of symmetry and of the fact that the underlying Markov chain cannot jump above $0$ (but with the advantage of having an elementary proof).

Another line of works considered persistence problems for additive functionals of \textit{continuous-time} Markov processes $(X_t)_{t\geq 0}$.
Let us mention~\cite{Goldman71} that considered the case of an integrated Brownian motion, and \cite{IK00,Prof21} who considered the $f$-integral of a Brownian motion or a skew-Bessel process respectively, for the (possibly asymmetric) functional $f(x) = |x|^{\gamma} (c_+\ind_{\{x>0\}} - c_- \ind_{\{x<0\}})$ for some $\gamma>-1$.
More recently, in~\cite{BBT23}, we have pushed further the existing techniques (based on a Wiener--Hopf decomposition of a bi-variate Lévy process associated to the problem): we obtained the sharp asymptotics for persistence probabilities for a wide class of Markov processes (including one-dimensional generalized diffusions, see~\cite{im96}) and of functions $f$.
In particular, \cite[Example~6]{BBT23} shows that the result applies to \textit{continuous-time} birth-death processes, giving the existence of the constant $\hat c_{\alpha}$ mentioned above.
The present article can be seen as a elementary approach to obtaining a sub-optimal result.

\section{Exchangeable and symmetric sequences: proof of Theorem~\ref{thm:exchang}}

Let us introduce some notation.
For $n\in \bbN$ and
for a fixed $x=(x_1, \ldots, x_n) \in \bbR^n$,
we define an exchangeable and symmetric vector
$\xi=(\xi_1,\ldots, \xi_n)$ with law denoted $\bP^{(x)}$ by
 permuting the coordinates of $(x_1, \ldots, x_n)$ by a random uniform permutation and by changing the signs of the coordinates uniformly at random.
More precisely, let $(\gep_i)_{1\leq i\leq n}$ be i.i.d.\ random variables with law $\bP(\gep_1 =1) =\bP(\gep_i=-1)=\frac12$
and let $\sigma$ be a random permutation of $\{1,\ldots, n\}$ with uniform distribution $\bP(\sigma =\nu) = \frac{1}{n!}$ for all $\nu\in \mathfrak{S}_n$,
independent of $(\gep_i)_{1\leq i \leq n}$:
we then define 
\begin{equation}
\label{def:X}
(\xi_1, \ldots, \xi_n) := (\gep_1 x_{\sigma(1)}, \ldots, \gep_n x_{\sigma(n)}) \,.
\end{equation}
Note that we can restrict to the case where $x_i \in \bbR_+$ for all $1\leq i\leq n$.

We then construct the random walk $S_k=\sum_{i=1}^k \xi_i$ for any $0\leq k \leq n$,
and we are interested in the persitence probabilities
\[
p_n(x) := \bP^{(x)} \big(S_k >0 \, \text{ for all } 1\leq k \leq n \big)
\quad \text{ and } \quad
\bar p_n(x) := \bP^{(x)} \big(S_k \geq 0 \, \text{ for all } 1\leq k \leq n \big) \,.
\]
We set by convention these probabilities equal to $1$ for $n=0$.
The following proposition is the essence of Theorem~\ref{thm:exchang}.

\begin{proposition}
\label{prop:exchangeable}
If $x=(x_1, \ldots, x_n)$ is such that 
\begin{equation}
\tag{H}
\label{hyp}
\sum_{i\in I} x_i \neq \sum_{j\in J} x_j 
\quad \text{ for all } I,J\subset \{1,\ldots, n\} \text{ with }I\neq J
\end{equation}
(by convention $\sum_{i\in \emptyset} x_i=0$), then we have that $p_n(x) = \bar p_n(x) = g(n):=\frac{1}{4^n} \binom{2n}{n}$; in particular it does not depend on $x$.
In general, we have $p_n(x) \leq g(n) \leq  \bar p_n(x)$.
\end{proposition}

Theorem~\ref{thm:exchang} follows directly from Proposition~\ref{prop:exchangeable}, noting that 
the law $\bP$ of $(\xi_1,\ldots, \xi_n)$, conditionally on $|\xi|=(|\xi_1|,\ldots, |\xi_n|)$, is $\bP^{(|\xi|)}$.
Let us also mention that a combinatorial proof of Proposition~\ref{prop:exchangeable} is given in Sections~2.3-2.4 of~\cite{Burns07} (see also \cite[Lem.~3.7]{BDGJS23}), in the spirit of that of~\cite{Sparre54}; our proof is quite simpler and does not use any combinatorial lemma.

%

\subsection{Proof of Proposition~\ref{prop:exchangeable}}
Our proof is greatly inspired by that of~\cite[Prop.~1.3]{DDG13}.
We start with the first statement:
we are going to prove by recurrence on $n$ that for any $x=(x_1,\ldots,x_n)$ that verifies the assumption~\eqref{hyp}, the quantities $p_n(x) = \bar p_n(x)$ do not depend on~$x$.
The statement is trivial for $n=1$ since we have $\bP^{(x)}(\xi_1>0)=\bP^{(x)}(\xi_1\geq 0) =\frac12$ if $x\neq 0$, so we directly proceed to the induction step.

Let us fix $n\geq 2$  and some $x=(x_1,\ldots,x_n)$ that verifies assumption~\eqref{hyp}.
We now apply a path decomposition used in~\cite{DDG13}.
Let $W = \min\{ k, S_k = \max_{1\leq i\leq n} S_i \}$: then, for any $\ell \in \{0,\ldots, n\}$, we have
\[
\begin{split}
\{W=\ell\} = \{ \xi_\ell>0, & \xi_\ell+\xi_{\ell-1}>0, \ldots, \xi_\ell+\cdots + \xi_1>0 \} \\
& \cap \{\xi_{\ell+1}\leq 0, \xi_{\ell+1}+\xi_{\ell+2}\leq 0, \ldots, \xi_{\ell+1}+\cdots +\xi_n \leq 0\} \,.
\end{split}
\]
Now, we have no independence at hand, but we can further decompose over permutations of $(x_1,\ldots, x_n)$
with a fixed image $I=\sigma(\{1,\ldots, \ell\})$: we get
\[
\begin{split}
&\bP^{(x)} (W =\ell)= \sum_{I\subset \{1, \ldots, n\},|I|=\ell}\ \sum_{\nu_1:\{1,\ldots,\ell\}\hookdoubleheadrightarrow I} \ \sum_{\nu_2:\{\ell+1,\ldots,n\}\hookdoubleheadrightarrow I^c}  \\[-3pt]
&\ \frac{1}{n!} \bP\Big( \sum_{i=0}^{j} \gep_{\ell-i} x_{\nu_1(\ell-i)} >0 \text{ for all } 0\leq j \leq \ell-1 \,,\, \sum_{i=1}^{j} \gep_{\ell+i} x_{\nu_2(\ell+i)} \leq 0 \text{ for all } 1\leq j\leq n-\ell \Big) \, \,,
\end{split}
\]
where $\nu:A\hookdoubleheadrightarrow B$ means that $\nu$ is a bijection from $A$ to $B$.
By independence and symmetry of the $(\gep_i)_{1\leq i \leq n}$, then recombining the sums over the permutations and using the exchangeability, we get that $\bP^{(x)} (W =\ell)$ is equal to
\[
 \sum_{I\subset \{1, \ldots, n\},|I|=\ell} \!\! \frac{\ell! (n-\ell)!}{n!} 
\bP^{(x_I)} \big( S_k >0 \, \text{ for all } 1\leq k \leq \ell  \big) \bP^{(x_{I^c})}\big(S_k \geq 0 \, \text{ for all } 1\leq k \leq n-\ell \big) \,,
\]
where we have defined $x_J = (x_i)_{i\in J}$ for any $J\subset \{1,\ldots, n\}$.
Since both $x_I$ and $x_{I^c}$ verify the assumption~\eqref{hyp}, 
we can apply the induction hypothesis:
we have that $p_{\ell}(x_I),p_{n-\ell}(x_{I^c})$ does not depend on $x_I, x_{I^c}$, for $1\leq \ell\leq n-1$.
Hence, denoting $p_\ell :=p_{\ell}(x_I)$, $p_{n-\ell} :=p_{n-\ell}(x_{I^c})$,  we end up with
$\bP^{(x)}(W =\ell) = p_{\ell} p_{n-\ell}$,
for any $\ell \in \{1, \ldots, n-1\}$.

Since we have $\bP^{(x)}(W=n) = \bP^{(x)}(\xi_n>0, \xi_n+\xi_{n-1}>0, \ldots, \xi_1>0)= p_n(x)$ by exchangeability and
$\bP^{(x)}(W=0) =\bP^{(x)}(S_k\leq 0 \text{ for all } 1\leq k \leq n) = \bar p_n(x)$ by symmetry,
we get
\[
1= \sum_{\ell=0}^n \bP^{(x)}(W=\ell) =  p_n(x)  + \sum_{\ell=1}^{n-1} p_\ell p_{n-\ell} + \bar p_n(x) \,.
\]
Using that $\sum_{i\in I} x_i \neq \sum_{j\in J} x_j$ for all $I,J\subset \{1,\ldots, n\}$ with $I\neq J$, we obtain that $p_n(x) =\bar p_n(x)$: the above identity shows that $p_n(x)=\bar p_n(x) =:p_n$ does not depend on $x$.

\smallskip
We can now determine the value of $p_n$, as done in~\cite{DDG13}. From the above, $(p_n)_{n\geq 0}$ satisfies the recursive relation $1=\sum_{\ell=0}^n p_{\ell} p_{n-\ell}$ for every $n\geq 0$ (recall that $p_0=1$). Constructing the generating function, we get that for any $|x|<1$,
\[
\sum_{n=0}^{\infty} x^n = \frac{1}{1-x} = \Big(\sum_{\ell=0}^{\infty} p_{\ell} x^\ell \Big)^2\,.
\]
Therefore, the generating function of $(p_{n})_{n\geq 0}$ is equal to $(1-x)^{-1/2}$, from which one deduces that $p_n = \frac{1}{4^n} \binom{2n}{n} = \frac{(2n-1)!!}{2n!!}$ for all $n\geq 0$.

\smallskip
For the general bounds, for any fixed $x$, for any $\delta >0$ fixed, one can choose $y=y(x,\delta)$ such that $x+y$ verifies the assumption~\eqref{hyp} and $\sum_{i=1}^n |y_i|\leq \delta$ (take e.g.\ a typical realization of i.i.d.\ random variables uniform in $[0,\delta/n]$).
Then, we clearly have that 
\[
 \bP^{(x+y)}\big( S_k \geq 0  \text{ for all } 1\leq k \leq n\big)
 \begin{array}{cl}
\geq  &\bP^{(x)}\big( S_k \geq \delta  \text{ for all } 1\leq k \leq n\big) \,,\\
\leq & \bP^{(x)}\big( S_k \geq -\delta  \text{ for all } 1\leq k \leq n\big) \,.
\end{array}
\]
Since $x+y$ verifies assumption~\eqref{hyp}, we get that the probability on the left-hand side does not depend on $x,y$ (and is equal to $\frac{(2n-1)!!}{2n!!}$), so 
\[
 \bP^{(x)}\big( S_k \geq \delta  \text{ for all } 1\leq k \leq n\big) \leq \frac{(2n-1)!!}{2n!!}\leq \bP^{(x)}\big( S_k \geq -\delta  \text{ for all } 1\leq k \leq n\big)\,.
\]
Since $\delta$ is arbitrary, letting $\delta \downarrow 0$ concludes the proof.\qed

\subsection{On the necessity of conditions~\ref{hyp:E}-\ref{hyp:S}}
\label{sec:counterexample}

Let us stress that in the construction~\eqref{def:X} of the exchangeable and symmetric vector $(\xi_1, \ldots, \xi_n)$,
the two assumptions are essential:

\smallskip
(a) \textit{The signs $(\gep_{i})_{1\leq i\leq n}$ need to be independent}. As a counter-example, take $(\gep_{i})_{1\leq i\leq n}$ uniform on  $\{\go \in \{-1,1\}^n, \sum_{i=1}^n \go_i \in \{2-n, n-2\} \}$, so that $\xi$ is still symmetric. 
Let $x =(x_1,\ldots, x_n)\in \bbR_+^n$ with $x_1\geq x_2 \geq \cdots \geq x_n$ (this is no restriction by definition of $\xi$). Then to have $S_k >0$ for all $1\leq k\leq n$, we need to have $\gep_1=+1$ (which happens with probability~$1/2n$) and then, since all other signs are $\gep_i =-1$, we need to place $x_1$ in the first position (which happens with probability~$1/n$): we get
\[
p_n(x)= \bP^{(x)} \big(S_k >0 \, \text{ for all } 1\leq k \leq n \big) = 
\begin{cases}
\frac{1}{2n^2} & \text{ if } x_1 >\sum_{i=2}^n x_i \,,\\
0 & \text{ if }x_1 \leq \sum_{i=2}^n x_i \,.
\end{cases}
\]
Therefore we obtain that $p_n(x)$ depends on $x$.

\smallskip
(b) \textit{The signs $(\gep_i)_{1\leq i \leq n}$ need to be independent from the permutation $\sigma$}. As a counter-example, take $n=3$ and $x=(x_1,x_2,x_3)\in (\bbR_+)^3$ with $x_1>x_2>x_3>0$: then the probability 
$\bP^{(x)} (S_1>0, S_2>0, S_3>0)$ is equal to (one needs to have $\gep_1=+1$)
\[
\begin{split}
&\bP(\gep_1=\gep_2=\gep_3 =1) + \bP(\gep_1=\gep_2=1, \gep_3=-1, x_{\sigma(3)}<x_{\sigma(1)}+x_{\sigma(2)}) \\
& + \bP(\gep_1=\gep_3=1, \gep_2=-1, x_{\sigma(2)}<x_{\sigma(1)}) + \bP(\gep_1=1, \gep_2=\gep_3=-1, x_{\sigma(2)} +x_{\sigma(3)}<x_{\sigma(1)}) \,.
\end{split}
\]
Now, if the joint distribution of $(\gep,\sigma)$
is such that $\gep =(\gep_i)_{1\leq i \leq 3}$ is uniform on $\{-1,1\}^3$ and $\bP( \sigma(3) = 3 \mid \gep_1=\gep_2=+1, \gep_3=-1)=1$, $\bP( \sigma(2) = 3 \mid \gep_1=\gep_2=+1, \gep_3=-1)=1$ and $\bP(\sigma(1)=1\mid \gep_1=+1, \gep_2=\gep_3=-1)=1$,
since $x_1>x_2>x_3>0$, we get that
\[
\bP^{(x)} (S_1>0, S_2>0, S_3>0)
= \frac14 + \frac18 \ind_{\{x_1> x_{2} +x_{3}\}} \,,
\]
which depends on $x$. It simply remains to see that the above conditions on the joint distribution of $(\gep,\sigma)$ can be satisfied, which can be checked by hand.

\section{Integrated birth and death chains}

\subsection{Persistence for integrated birth-death chains: Proof of Theorem~\ref{thm:integrated}}

Let us define iteratively $\tau_0=0$ and, for $k\geq 1$, $\tau_k = \min\{n > \tau_{k-1}, X_n=0\}$.
Then, for $k\geq 1$, we define the random variable
\[
\xi_k = \sum_{i=\tau_{k-1}+1}^{\tau_k} f(X_i) \,,
\]
\textit{i.e.}\ the contribution of the $k$-th excursion of $(X_n)_{n\geq 0}$ to the $f$-integrated Markov chain.
Note that by the Markov property, the $(\xi_k)_{k\geq 1}$ are i.i.d.
We can therefore write $\zeta_n := \sum_{i=1}^n f(X_i)$ as
\[
\zeta_n = \sum_{k=1}^{L_n} \xi_k + W_n, \qquad W_n = \sum_{i=\tau_{L_n}+1}^{n} f(X_i)
\]
where we recall that $L_n =\sum_{i=1}^n\ind_{\{X_i=0\}}$ is the local time at $0$.
Since there is no change of sign during an excursion (recall that $(X_n)_{n\geq 0}$ is a birth-death chain), 
we have that $(\zeta)_{n\geq 0}$ is monotonous on each interval $(\tau_{k-1}, \tau_k]$.
We therefore get the following bounds.

Removing the positivity condition on the last segment $(\tau_{L_n}, n]$, we get
\begin{equation}
\label{eq:upper}
\bP(\zeta_k > 0 \text{ for all }1\leq k\leq n )
\leq  \bP\Big( \sum_{k=1}^{\ell} \xi_k  > 0 \text{ for all } 1\leq \ell \leq L_n\Big) \,.
\end{equation}
On the other hand, imposing that the excursion straddling over $n$ is non-negative (with probability $1-p_0$) so that $W_n\geq 0$, and using that the sign of $W_n$ is independent from the past by symmetry, we get
\begin{equation}
\label{eq:lower}
\bP(\zeta_k >  0 \text{ for all }1\leq k\leq n )
\geq (1-p_0) \bP\Big(\sum_{k=1}^{\ell} \xi_k > 0 \text{ for all } 1\leq \ell \leq L_n\,, L_n\geq 1\Big) \,.
\end{equation}

The difficulty now to study $\bP(\sum_{k=1}^{\ell} \xi_k > 0 \text{ for all } 1\leq \ell \leq L_n)$ is that the number of terms is random and that $\xi_i$ and $L_n$ are not independent. 
However, for any $0\leq m\leq j\leq n$,  conditionally on $\{L_n = m, \tau_{m} = j\}$, the random variables $(\xi_i)_{1\leq i \leq m}$, though not independent, are easily seen to be exchangeable and sign-invariant (thanks to the Markov property and the fact that the chain $(X_n)_{n\geq 0}$ is symmetric).
We can therefore apply Theorem~\ref{thm:exchang} and Remark~\ref{rem:bounds1} to obtain that
\[
p_0 g(m-1) \leq \bP\Big( \sum_{k=1}^{\ell} \xi_k  > 0 \text{ for all } 0\leq \ell \leq m \, \Big|\,  L_n = m, \tau_m = j\Big) \leq  g(m)\,.
\]
We therefore get that, conditioning on $L_n,\tau_{L_n}$,
\begin{align*}
p_0 g(L_n-1) \leq \bP\Big( \sum_{k=1}^{\ell} \xi_k  > 0 \text{ for all } 0\leq \ell \leq L_n \;\Big| \; L_n, \tau_{L_n}\Big)
 \leq g(L_n) \,,
\end{align*}
which gives the desired bounds thanks to~\eqref{eq:upper}-\eqref{eq:lower}, by taking the expectation.

As far as bridges are concerned, we have the identity
\begin{equation}
\label{eq:equal}
\bP\big(\zeta_k > 0 \text{ for all }1\leq k\leq n, X_n=0 \big)
=  \bP\Big( \sum_{k=1}^{\ell} \xi_k  > 0 \text{ for all } 1\leq \ell \leq L_n , X_n=0\Big) \,,
\end{equation}
so similarly as above we end up with the following conditional expectation
\[
p_0 g(L_n-1) \ind_{\{X_n=0\}} \leq \bP\Big( \sum_{k=1}^{\ell} \xi_k  > 0 \text{ for all } 0\leq \ell \leq L_n , X_n=0\;\Big| \; L_n, \tau_{L_n}\Big)
\leq g(L_n) \ind_{\{X_n=0\}}\,,
\]
which gives the desired bound, again by taking the expectation.\qed



\begin{remark}
We also have, with the same line of proof and using Theorem~\ref{thm:exchang}  
\[
\begin{split}
\bP\big(\zeta_k \geq 0 \text{ for all } 1\leq k \leq n \big) & \geq (1-p_0) \bE\big[ g(L_n)\big] \,,\\
\bP\big(\zeta_k \geq 0 \text{ for all } 1\leq k \leq n , X_n =0\big)& \geq (1-p_0) \bE\big[ g(L_n)\, \ind_{\{X_{n}=0\}} \big] \,.
\end{split}
\]
\end{remark}

\subsection{Proof of Proposition~\ref{lem:Ln}}

First of all, let us notice that since $ (\pi(m+1/2))^{-1/2} \leq g(m) \leq (\pi m)^{-1/2}$  for all $m\geq 1$, we only have to obtain the asymptotics of 
\begin{equation}
\label{eq:expectLn}
\bE\big[ (L_n)^{-1/2} \ind_{\{L_n \geq 1\}} \big] 
\qquad \text{ and } \qquad
\bE\big[ (L_n)^{-1/2} \ind_{\{X_n=0\}}\ind_{\{L_n \geq 1\}} \big] \,.
\end{equation}
The same asymptotics with $(L_n)^{-1/2}$ replaced with $(L_n+1/2)^{-1/2}$ follow analogously.

\medskip
\noindent
{\it Step 1. Convergence of $b_n^{-1}L_n$.} Let us show that under the assumptions of Proposition~\ref{lem:Ln}, we have that $(b_n^{-1} L_n)_{n\geq 0}$ converges in distribution to a random variable $\mathcal{X}$.
This is obvious in the positive recurrent case with $b_n=n$, thanks to the ergodic theorem, with $\mathcal{X} = \bE[\tau_1]^{-1}$.
In the null-recurrent case, we use the fact that 
$\bP( b_n^{-1}L_n < t) = \bP(\tau_{\lceil t b_n \rceil} >n)$ to deduce the convergence in law of $(b_n^{-1} L_n)_{n\geq 0}$ from that of $(\tau_k)_{k\geq 1}$.
Let $(a_n)_{n\ge 1}$ be a sequence such that  $\bP(\tau_1 > a_n)\sim n^{-1}$ as $n\to\infty$. Then we consider the three cases $\alpha=1$, $\alpha\in (0,1)$ and $\alpha=0$ separately.


\smallskip
\noindent
\textit{Case $\alpha=1$}. We then have that $a_n^{-1}(\tau_n- n \mu(a_n))$ converges in distribution to a $1$ stable law, see~\cite[IX.8]{FellerII} (see Eq.~(8.15) for the centering): in particular, $\frac{\tau_n}{n\mu(a_n)} $ converges in probability to~$1$ (recall that $\mu(n):=\bE[\tau_1 \ind_{\{\tau_1 \leq n\}}]$, which is slowly varying).
Setting $b_n$ such that $b_n \mu(a_{b_n}) \sim n$, then we have that $\frac{\tau_{tb_n}}{n}$ converges to $t$ in probability, so $b_n^{-1} L_n$ converges in probability to~$1$.
Then, we simply have to notice that if $b_n$ is given by the above relation then we  $\mu(a_{b_n})\sim \mu(n)$, see e.g.~\cite[Lem.~4.3]{B18}, so $b_n\sim n/\mu(n)$ as defined in~\eqref{eq:persist2}.

\smallskip
\noindent
\textit{Case $\alpha \in (0,1)$.}  We then have that $\frac{1}{a_n} \tau_n$ converges in distribution to $\kappa_{\alpha} \mathcal{Z}_{\alpha}$, with $\kappa_{\alpha}= \Gamma(1-\alpha)^{1/\alpha}$ and $\mathcal{Z}_{\alpha}$ a one-sided $\alpha$-stable random variable~$\mathcal Z$ with Laplace transform $e^{-t^\alpha}$, see~\cite[XIII.6]{FellerII}.
Then setting $b_n:= n^{\alpha}/\ell(n)$ as in~\eqref{eq:persist2}, we get that $a_{b_n}\sim n$ as $n\to\infty$; indeed, $\bP(\tau_1 > n)\sim b_n^{-1}\sim\bP(\tau_1 > a_{b_n})$ by definition of $a_n$.
Using that $a_n$ is regularly varying with exponent $1/\alpha$, we get that $a_{t b_n} \sim t^{1/\alpha} a_{b_n} \sim t^{1/\alpha} n$: we therefore get that $\lim_{n\to\infty} \bP(\tau_{tb_n} >n) = \bP( \kappa_{\alpha}\mathcal Z_{\alpha} > t^{-1/\alpha})$. 
This shows that $b_n^{-1}L_n$ converges in distribution to $(\kappa_{\alpha} \mathcal{Z}_{\alpha})^{-\alpha}$.

\smallskip
\noindent
\textit{Case $\alpha=0$.} We then have that $n \ell(\tau_n)$ converges in distribution to an $\mathrm{Exp}(1)$ variable, see \cite[Thm.~8]{D52}.
Let $b_n = 1/\ell(n)$ as in~\eqref{eq:persist2} and assume (without loss of generality) that $\ell$ is strictly decreasing.
Then, for any $t>0$,
\[
\bP( b_n^{-1} L_n < t) = \bP(\tau_{tb_n} >n) = \bP\big( t b_n \ell(\tau_{tb_n}) < t b_n \ell(n) \big) \xrightarrow[n\to\infty]{}  1-e^{-t} \,,
\]
which shows that $b_n^{-1}L_n$ converges in distribution to an $\mathrm{Exp}(1)$ variable.

\medskip
\noindent
\textit{Step 2a. Convergence of the expectation $\bE[ (L_n)^{-1/2} \ind_{\{L_n\geq 1\}}]$.}
We have that
\begin{align}
 \bE\big[ (L_n)^{-1/2} \ind_{\{L_n\geq 1\}} \big] &  = \int_0^1 \bP\big( (L_n)^{-1/2} > t, \: L_n \geq 1 \big) \dd t \nonumber \\
  & = b_n^{-1/2} \int_{b_n^{-1}}^{\infty} \bP\big(b_n^{-1}\leq  b_n^{-1} L_n <  u \big)   \tfrac12 u^{-3/2}  \dd u \,. \label{eq:integral}
\end{align}
Since by Step~1 we have that $\bP(b_n^{-1}\leq b_n^{-1} L_n \leq  u )$ converges to $\bP(\mathcal{X}\leq u)$ as $n\to\infty$, we simply need some uniform bound on $\bP( b_n^{-1} L_n \leq  u )$ in order to be able to apply dominated convergence; note that we only need a bound for $u\leq 1$, otherwise we simply bound the probability by~$1$.

\smallskip
\noindent
\textit{(i) Positive recurrent case.}
In that case $b_n=n$, so we easily get that uniformly in $n\geq 1$
\[
\bP( b_n^{-1} L_n \leq  u ) = \bP(\tau_{u n} >n ) \leq n^{-1} \bE[\tau_{u n}] = u \bE[\tau_1]\,,
\]
by Markov's inequality.
Since $u \times u^{-3/2}$ is integrable on $(0,1]$, we can apply dominated convergence to get that
\[
\lim_{n\to\infty} \int_{b_n^{-1}}^{\infty} \bP\big(b_n^{-1}\leq  b_n^{-1} L_n <  u \big)   \tfrac12 u^{-3/2}  \dd u  = \int_{0}^{\infty} \bP\big( \mathcal{X} \leq  u \big)   \tfrac12 u^{-3/2}  \dd u  = \bE[\mathcal{X}^{-1/2}] = \bE[\tau_1]^{1/2} \,,
\]
where we have used that $\mathcal{X}=\bE[\tau_1]^{-1}$ in the positive recurrent case.
Combined with~\eqref{eq:integral}, this gives the desired asymptotics.

\medskip
\noindent
\textit{(ii) Null-recurrent case.}
In the null-recurrent case, we use that 
\[
\bP( b_n^{-1} L_n \leq  u ) = \bP(\tau_{u b_n} >n )
\leq  C \, u b_n \,\bP(\tau_1 >n) \,,
\]
for some constant $C>0$.
The last bound is standard and falls in the ``big-jump'' phenomenon (note that the event $\tau_{t b_n} > n$ is an upper large deviation).
If $\alpha=1$, this is given by~\cite[Thm.~1.2]{Nag79}; note also that in this case $\lim_{n\to\infty} b_n \bP(\tau_1 >n) =0$, thanks to~\cite[Prop.~1.5.9.a.]{BGT89}.
If $\alpha\in (0,1)$, this is contained in~\cite[Thm.~1.1]{Nag79}; note that by definition of $b_n$ we have that $b_n \bP(\tau_1 >n)$ remains bounded in this case.
In the case $\alpha=0$, this is due to~\cite{NV08} (see also~\cite{AB16a}); and $b_n \bP(\tau_1 >n)$ also remains bounded in this case, by definition of $b_n$.
This shows that in all cases we have $\bP( b_n^{-1} L_n \leq  u ) \leq C u$, so as above we can apply dominated convergence to get that
\[
\lim_{n\to\infty} \int_{b_n^{-1}}^{\infty} \bP\big(b_n^{-1}\leq  b_n^{-1} L_n <  u \big)
\tfrac12 u^{-3/2}  \dd u  = \int_{0}^{\infty} \bP\big( \mathcal{X} \leq  u \big)   \tfrac12 u^{-3/2}  \dd u  = \bE[\mathcal{X}^{-1/2}]\,,
\]
with $\mathcal{X}=1$ if $\alpha=1$, $\mathcal{X}=(\kappa_{\alpha} \mathcal{Z})^{-\alpha}$ if $\alpha\in (0,1)$ and $\mathcal{X} \sim \mathrm{Exp}(1)$ if $\alpha=0$.
Combined with~\eqref{eq:integral}, this gives the desired asymptotics.

\medskip
\noindent
\textit{Step 2b. Convergence of the expectation $\bE[ (L_n)^{-1/2} \ind_{\{X_n=0\}}\ind_{\{L_n\geq 1\}}]$.}
In the positive recurrent case, for any $\varepsilon>0$, let us set $A_{\gep, n} = \{n^{-1}L_n \in (1-\gep,1+\gep)\bE[\tau_1]^{-1}\}$.
 With the same argument as above, we have that 
\[
\lim_{n\to\infty} \sqrt{n}\bE[ (L_n)^{-1/2} \ind_{\{X_n=0\}}\ind_{\{L_n\geq 1\}}\ind_{ A_{\gep,n}^c}] = 0 \,.
\]
Indeed, we first observe that, as in \eqref{eq:integral}, we have
\[
 \sqrt{n}\bE\big[ (L_n)^{-1/2} \ind_{\{L_n\geq 1\}} \ind_{A_{\gep, n}^c}\big] = \int_{1/n}^\infty \bP\big(n^{-1}\leq  n^{-1} L_n <  u, \: A_{\gep, n}^c \big) \tfrac12 u^{-3/2}  \dd u.
\]
Since for every $\epsilon > 0$, $\bP(A_{\gep, n}^c)$ vanishes as $n\to\infty$, the above limit holds using the same dominated convergence argument as in the previous step. Therefore, we only need to consider the remaining term $\bE[ (L_n)^{-1/2} \ind_{\{X_n=0\}}\ind_{A_{\gep,n}}]$: by definition of $A_{\gep, n}$, we have that
\[
(1+\gep)^{-1/2} \bP(X_n=0, A_{n,\gep}) \leq  \sqrt{n}\bE[ (L_n)^{-1/2} \ind_{\{X_n=0\}}\ind_{A_{\gep,n}^c}]  \leq  (1-\gep)^{-1/2} \bP(X_n=0) \,.
\]
Now, we have $\lim_{n\to\infty} \bP(X_n=0) = 1/\bE[\tau_1]$ by the convergence to the stationary distribution, and we also have $\lim_{n\to\infty} \bP(X_n=0, A_{\gep, n}) = 1/\bE[\tau_1]$
since $\lim_{n\to\infty}\bP(X_n=0, A_{\gep, n}^c) =0$.
This concludes the proof of~\eqref{eq:persist1}.

\smallskip
For the null-recurrent case, we focus on the case $\alpha\in (0,1)$. We write
\begin{equation}
\bE\big[ (L_n)^{-1/2} \ind_{\{X_n=0\}}\ind_{\{L_n\geq 1\}}\big] = \sum_{k=1}^{\infty} k^{-1/2} \bP(L_n=k, X_n=0) = \sum_{k=1}^{\infty} k^{-1/2} \bP(\tau_k=n) \,.
\end{equation}
Our proof is quite standard and follows the line of~\cite[\S3]{D97} (with some adaptation).
We fix $\gep>0$ and we decompose the sum into three parts:
\[
P_1 := \sum_{1\leq k < \gep b_n} k^{-1/2} \bP(\tau_k=n) \,,
\quad
P_2 := \sum_{k= \gep b_n}^{\gep^{-1} b_n} k^{-1/2} \bP(\tau_k=n) \,,
\quad
P_3 := \sum_{k> \gep^{-1} b_n} k^{-1/2} \bP(\tau_k=n)\,.
\]
The main contribution comes from the term $P_2$. Gnedenko's local limit theorem~\cite{GK54} gives that $\delta_k :=\sup_{x\in \mathbb Z} \big| a_k \bP(\tau_k =x) - g_{\alpha}\big( \frac{x}{a_k} \big)  \big|$ goes to $0$ as $k\to\infty$, where $g_{\alpha}$ is the density of the limiting $\alpha$-stable distribution $\kappa_{\alpha} \mathcal{Z}_{\alpha}$.
We therefore have that
\[
\bigg| P_2 - \sum_{k= \gep b_n}^{\gep^{-1} b_n}  \frac{k^{-1/2} }{a_k} g_{\alpha}\Big( \frac{n}{a_k} \Big) \bigg|  \leq  \bar \delta_{\gep, n} \sum_{k= \gep b_n}^{\gep^{-1} b_n} k^{-1/2}  \frac{1}{a_k} \,,
\]
with $\bar \delta_{\gep, n} = \sup_{k\geq \gep b_n}  \delta_k$ going to $0$ as $n\to\infty$.
Since $k^{-1/2}/a_k$ is regularly varying, we have that $k^{-1/2}/a_k \leq C_{\gep} b_n^{-1/2}/a_{b_n}$ uniformly in $k\in [\gep b_n, \gep^{-1}b_n]$, so using also that $a_{b_n}\sim n$, we get 
\[
\bigg| P_2 - \sum_{k= \gep b_n}^{\gep^{-1} b_n}  \frac{k^{-1/2} }{a_k} g_{\alpha}\Big( \frac{n}{a_k} \Big) \bigg|  \leq  \bar \delta_{\gep n} \gep^{-1} b_n \times C'_{\gep} b_n^{-1/2}/n = o(b_n^{1/2}/n) \,.
\]
We now focus on the remaining sum.
Since $a_k$ is regularly varying with exponent $1/\alpha$ we get that $a_k = (1+o(1))  (k/b_n)^{1/\alpha} a_{b_n}$ as $n\to\infty$ with the $o(1)$ uniform in $k\in [\gep b_n, \gep^{-1}b_n]$ (it may depend on $\gep$).
Using also that $a_{b_n}\sim n$ we get that $\frac{n}{a_k} = (1+o(1)) (k/b_n)^{-1/\alpha}$:
thanks to the uniform continuity of $g_{\alpha}$ on $[\gep, \gep^{-1}]$, we get that
\[
\sum_{k= \gep b_n}^{\gep^{-1} b_n}  \frac{k^{-1/2} }{a_k} g_{\alpha}\Big( \frac{n}{a_k} \Big) = (1+o(1)) \sum_{k= \gep b_n}^{\gep^{-1} b_n} k^{-1/2} \frac1n \Big( \frac{k}{b_n}\Big)^{-\frac{1}{\alpha}}  g_{\alpha}\Big(  \Big( \frac{k}{b_n}\Big)^{-1/\alpha} \Big) \,.
\]
Therefore, a Riemann sum approximation for the last identity yields that
\[
\lim_{n\to\infty} n b_n^{-1/2} P_2 = \int_{\gep}^{\gep^{-1}} t^{-\frac{1}{2}-\frac{1}{\alpha}} g_{\alpha}\big(t^{-1/\alpha}\big) \dd t \,.
\]

One can also control $P_3$ thanks to the local limit theorem: there is a constant $C$ such that $\bP(\tau_k=n) \leq C/a_k$ for all $k\geq 1$, so that using the same Riemann sum approximation as above,
\[
P_3 \leq C \sum_{k> \gep^{-1} b_n} \frac{k^{-1/2}}{a_k}  \leq C' \frac{ b_n^{1/2}}{n} \int_{\gep^{-1}}^{\infty} t^{-\frac{1}{2}-\frac{1}{\alpha}} \dd t = C_{\alpha}'\, \gep^{\frac{1}{\alpha}-\frac12} \frac{ b_n^{1/2}}{n}  \,,
\]
where we have used that $\frac12+\frac1\alpha >1$ in the last identity (since $\alpha\in(0,1)$).
This shows that $\limsup_{n\to\infty} n b_n^{-1/2} P_3$ goes to $0$ as $\gep\downarrow 0$.

It remains to deal with $P_1$: we need to show that, as for $P_3$, 
\begin{equation}
\label{eq:termP1}
\lim_{\gep\downarrow 0}\limsup_{n\to\infty}  nb_n^{-1/2} P_1 =0 \,.
\end{equation}
With this at hand, letting $n\to\infty$ then $\gep\downarrow 0$, we would then get that 
\[
\lim_{n\to\infty}  nb_n^{-1/2} \sum_{k=1}^{\infty} k^{-1/2} \bP(\tau_k=n) = \int_{0}^{\infty}  t^{-\frac{1}{2}-\frac{1}{\alpha}} g_{\alpha}\big(t^{-1/\alpha}\big) \dd t  = \alpha \bE\big[ \big(\kappa_{\alpha} \mathcal{Z}_{\alpha} \big)^{-\alpha/2} \big] \,,
\]
the last identity following from a simple change of variable. Since $b_n = n^{\alpha} \ell(n)^{-1}$, this shows~\eqref{eq:persistlocal}.

\smallskip
Now, it only remains to prove~\eqref{eq:termP1}, and this is where the condition $\alpha>2/3$ will appear. We need to use local large deviations to bound $\bP(\tau_k=n)$.
First, a general bound is given in~\cite[Thm.~2.3]{B19} (see also \cite[Thm.~6.1]{CD19}): wthere is a constant $C>0$ such that $\bP(\tau_k =n )\leq \frac{C}{a_k} k \ell(n) n^{-\alpha}$ for all $n\geq k$. This gives that, for $\alpha >2/3$,
\[
P_1 =  \sum_{1\leq k < \gep b_n} k^{-1/2} \bP(\tau_k=n)\leq C n^{-\alpha} \ell(n) \sum_{k=1}^{\gep b_n} \frac{k^{1/2}}{a_k} \leq C' b_n^{-1} \times   \frac{ (\gep b_n)^{3/2}}{a_{\gep b_n}} \,.
\]
The fact that $\alpha>2/3$ is crucial in the last inequality, where we have used that $k^{1/2}/a_k$ is regularly varying with exponent $\frac{1}{2} -\frac{1}{\alpha} >-1$.
Now, since $a_{\gep b_n} \sim \gep^{1/\alpha} a_{b_n} \sim \gep^{1/\alpha} n$, we get that
$P_3 \leq C'' \gep^{\frac{3}{2} - \frac{1}{\alpha}} n^{-1} b_n^{1/2}$ for all $n$. This shows that $\limsup_{n\to\infty} nb_n^{-1.2}P_1 \leq C'' \gep^{\frac{3}{2} - \frac{1}{\alpha}}$, so that~\eqref{eq:termP1} holds (again, using that $\frac{3}{2} - \frac{1}{\alpha}>0$ for $\alpha>2/3$).

In the case where $\alpha\in (0,2/3]$, one needs an extra assumption to obtain a better local large deviations estimate. 
From~\cite[Thm.~2]{D97} (or \cite[Thm.~2.4]{B19}), the condition $\bP(\tau_1 =n)\leq C n^{-(1+\alpha)} \ell(n)$ ensures that there is a constant $C'>0$ such that $\bP(\tau_k =n )\leq C' k \ell(n) n^{-(1+\alpha)}$.
Hence, in this case, one ends up with
\[
P_1= \sum_{1\leq k < \gep b_n} k^{-1/2} \bP(\tau_k=n) \leq C' n^{-(1+\alpha)} \ell(n) \sum_{k=1}^{\gep b_n} k^{1/2} \leq C'' n^{-1} b_n^{-1} \times  (\gep b_n)^{3/2} = C''  \gep^{3/2} n^{-1} b_n^{1/2}.
\]
This shows that $\limsup_{n\to\infty} nb_n^{-1.2}P_1 \leq C'' \gep^{3/2}$ so that~\eqref{eq:termP1} holds.
This concludes the proof of~\eqref{eq:persistlocal}.
\qed

\paragraph*{Acknowledgements.} We are grateful for the numerous discussions on this subject with our colleagues from LPSM; we would like to thank in particular Thomas Duquesne, Nicolas Fournier and Camille Tardif.
We also want to thank Igor Kortchemski for pointing out reference~\cite{BDGJS23}.

{
\bibliographystyle{abbrv}
\bibliography{biblio.bib}
}

{
  \bigskip
  \footnotesize

  Q.~Berger, \textsc{Sorbonne Universit\'e, Laboratoire de Probabilit\'es, Statistique et Modélisation, 75005 Paris, France} and 
  \textsc{DMA, École Normale Supérieure, Université PSL, 75005 Paris, France.}\par\nopagebreak
  \textit{E-mail address}: \texttt{quentin.berger@sorbonne-universite.fr}

  \medskip

  L.~B\'ethencourt, \textsc{Sorbonne Universit\'e, Laboratoire de Probabilit\'es, Statistique et Modélisation, 75005 Paris, France.}\par\nopagebreak
  \textit{E-mail address}: \texttt{loic.bethencourt@sorbonne-universite.fr}
}

\end{document}